\documentclass[12pt]{amsart}
\usepackage{amssymb}
\usepackage{eucal}
\theoremstyle{plain}
\newtheorem{theorem}{Theorem}
\newtheorem{corollary}{Corollary}

\newtheorem{lemma}{Lemma}

\theoremstyle{definition}

\newtheorem{remark}{Remark}
\theoremstyle{remark}

\numberwithin{equation}{section}
\newcommand{\p}{\partial}

\newcommand{\Rn}{\mathbb R^n}

\begin{document}
\title[Convexity of free solutions of Schr\"odinger equation with\dots]{Convexity properties of solutions to the free Schr\"odinger equation with Gaussian decay}
\author{L. Escauriaza}
\address[L. Escauriaza]{Universidad del Pa{\'\i}s Vasco / Euskal Herriko
Unibertsitatea\\Dpto. de Matem\'aticas\\Apto. 644, 48080 Bilbao, Spain.} 
\email{mtpeszul@lg.ehu.es}
\thanks{The first and fourth authors are supported  by MEC grant MTM2007-62186. The second and third authors were supported by NSF grants.}
\author{C.E. Kenig}
\address[C.E. Kenig]{Department of Mathematics\\University of Chicago\\5734 S. University Avenue\\Chicago,
Illinois 60637, USA.} 
\email{cek@math.uchicago.edu}
\author{G. Ponce}
\address[G. Ponce]{Department of Mathematics\\ South Hall, Room 6607\\ University of California, Santa Barbara\\CA 93106, USA.} 
\email{ponce@math.ucsb.edu}
\author{L. Vega}
\address[L. Vega]{Universidad del Pa{\'\i}s Vasco / Euskal Herriko
Unibertsitatea\\Dpto. De matem\'aticas\\Apto. 644, 48080 Bilbao, Spain.} 
\email{mtpvegol@lg.ehu.es}
\begin{abstract} We study convexity properties of solutions to the free Schr\"odinger equation with Gaussian decay.
\end{abstract}
\maketitle

\section{Introduction}\label{S:1}

The main purpose of this note is to study convexity properties of two classes of functions.
The first one is the class of functions which, together with their Fourier transform, have Gaussian
decay. The second class is the class of functions for which the free Schr\"odinger evolution has Gaussian decay at two different times. (We will see later that, in fact, the two classes coincide). As motivation for the study of these issues, we recall the well known \lq\lq uncertainty principle" due to G. H. Hardy (see \cite{StSh}) :

(A) {\it Assume $n=1$, $f(x) = O(e^{-Ax^2})$ and $\hat f(\xi)=O(e^{-B\xi^2})$. If \newline
$A B>1/4$, then $f\equiv 0$. Moreover, if $A B=1/4$, then  $f(x)=c e^{-Ax^2}$, for some constant $c$.  In \cite{SiSu} this result was extended  to higher dimensions, with $x^2$ and $\xi^2$ replaced by $|x|^2$ and $|\xi|^2$}.

Here \[\widehat f(\xi)=(2\pi)^{-\frac n2}\int_{\Rn}e^{-ix\cdot\xi}f(x)\,dx.\]

As pointed out in \cite{EsKePoVe}, this result has an equivalent formulation for the free Schr\"odinger
equation. Thus, consider the IVP for the free Schr\"odinger equation 
\begin{equation*}
\begin{cases}
\begin{aligned}
\label{IVP}
&\partial_tu= i\Delta u,
\;\;\;x\in\mathbb R^n,\;t\in\mathbb R,\\
&u(x,0)=u_0(x),
\end{aligned}
\end{cases}
\end{equation*}
whose solution $u(x,t)=e^{it\Delta}u_0(x)$ can be written as
\begin{equation*}
\begin{aligned}\label{de-da}
&u(x,t)=  (e^{-i|\xi|^2t} \widehat u_0)^\lor(x)=\int_{\mathbb R^n} \frac{e^{i|x-y|^2/4t}}{(4\pi i t)^{n/2}}\, u_0(y)\,dy\\
=&\frac{e^{i|x|^2/4t}}{(4\pi i t)^{n/2}} \int_{\mathbb R^n}e^{-2ix\cdot \xi/4t} e^{i|y|^2/4t} u_0(y)\,dy
= \frac{e^{i|x|^2/4t}}{(2 i t)^{n/2}}\;
\widehat{\;(e^{i|\cdot|^2/4t}u_0)\,}\left(\frac{x}{2 t}\right),
\end{aligned}
\end{equation*}
i.e. if $c_t=(2 i t)^{n/2}$, then
\begin{equation}
\label{dada}
c_t e^{-i|x|^2/4t} \,u(x,t) = (e^{i|\cdot|^2/4t}u_0)^{\land}\left(\frac{x}{2 t}\right).
\end{equation}

 Formula \eqref{dada} tells us that  $e^{-i|x|^2/4t} \,u(x,t)$ is a multiple of a rescaled Fourier transform
 of $e^{i|y|^2/4t}u_0(y)$, thus Hardy's result can be restated in terms of the Schr\"odinger equation as:
 
 (A$'$) {\it If $u_0(x)=O(e^{-|x|^2/\beta^2})$, $u(x,t)=e^{it\Delta}u_0(x)=O(e^{-|x|^2/\alpha^2})$, $t>0$ and $\alpha \beta<4t$, then $u\equiv 0$ in $\mathbb R^n\times[0,t]$. Moreover, if $\alpha \beta = 4t$, then $u$ is the solution with initial data, $\omega e^{-\left(\frac{1}{\beta^2}+\frac{i}{4t}\right)|x|^2}$, for some complex number $\omega$.} 
 \vskip.05in 
 We will also recall the following extension of (A) established in \cite{SiSu}: 
  \vskip.05in
 (B) {\it If $e^{A_1|\,\cdot\,|^2}f\in L^p(\mathbb R)$ and $e^{A_2|\,\cdot\,|^2}\widehat f\in L^q(\mathbb R)$, $p,\,q\in[1,\infty]$, with at least one of them finite and with $A_1 A_2\geq  1/4$, then $f\equiv 0$,} 
 \vskip.05in 
\noindent and the Beurling-H\"ormander result in \cite{Ho}:  
 \vskip.05in 
(C) {\it If $f\in L^1(\mathbb R)$ and $\int_{\mathbb R} \int_{\mathbb R} |f(x)| |\widehat f(\xi)| e^{ |x\,\xi|}\,dx\,d\xi<\infty$, then $f\equiv 0$.} 
\vskip.05in 
Because of \eqref{dada}, (B) and (C) can be rephrased as:
\vskip.05in
(B$'$) {\it If $u_0\in L^p(e^{p\, x^2/\beta^2}\,dx)$ and $e^{it\partial_x^2}u_0\in L^q(e^{q \,x^2/\alpha^2}\,dx)$, $p,\,q\in[1,\infty]$, with at least one of them finite and $4t\geq \alpha \beta$, then $u_0\equiv 0$,}
\vskip.05in
\noindent
and
\vskip.05in 
(C$'$) {\it If $u_0\in L^1(\mathbb R)$ and $\int_{\mathbb R} \int_{\mathbb R} |u_0(x)| | e^{it\partial_x^2}u_0(y)| e^{ |x y|/2t}\,dx \,dy<\infty$, then $u_0\equiv 0$.}
\vskip.1in

  In \cite{EsKePoVe}, we obtained the following results which can be seen as variants of Hardy's uncertainty principle, in the context of Schr\"odinger equations with lower order variable coefficients and for non-linear Schr\"odinger equations.

 
(I) \it Let $u\in C([0,1]:H^2(\mathbb  R^n))$ be a strong solution of 
 \begin{equation}
 \label{1.1c}
 i\partial_t u + \Delta u =V u,
\end{equation}
with 
$$
V:\mathbb  R^n\times [0,1]\to\mathbb  C,\;\;V\in L^{\infty}(\mathbb  R^n_x
 \times[0,1]),\;\;\;\nabla_x V\in L^1_t([0,1]:L^{\infty}(\mathbb  R^n)),
 $$  
and
$$
\lim_{r\uparrow \infty} \,\| V\|_{L^1_tL^{\infty}_x\{|x|>r\}}
= \lim_{r\uparrow \infty}\,\int_0^1\,\sup_{|x|>r}\,|V(x,t)|\,dt=0.
 $$
 There exists  $
 c_0=c_0(n; \|u\|_{L^{\infty}_tH^2_x};\|V\|_{L_{t,x}^{\infty}};\|\nabla_xV\|_{L^1_tL^{\infty}_x})>0 $
 such that if 
 \begin{equation}
\label{1.3}
 u_0=u(0),\;\;\;u_1=u(1)\in H^1(e^{a|x|^{2}}dx),
 \end{equation}
 with $a\geq c_0$, 
 then $u\equiv 0$. \rm
 
 
  \vskip.15in
 
 (II) \it Let $u_1,\,u_2\in C([0,1]:H^k(\mathbb  R^n))$, $\,k>n/2+1$ be solutions of
 \begin{equation}
 \label{1.2c}
 i\p_t u + \Delta u + F(u,\overline u)=0,
\end{equation}
with 
$$
F:\mathbb  C^2\to \mathbb  C,\;\;F\in C^k,\;\;F(0)=\partial_uF(0)=\partial_{\bar u}F(0)=0.
$$
There exists  $ c_0=c_0(n; \|u_1\|_{L^{\infty}_tH^2_x};\|u_2\|_{L^{\infty}_tH^2_x};\|F\|_{C^k})>0
 $
 such that if  
 \begin{equation}
\label{1.7}
  u_1(0)-u_2(0),\;\;\;\;\;\;\;\;u_1(1)-u_2(1)\in H^1(e^{a|x|^{2}}dx),
 \end{equation}
 with $a\geq c_0$,
 then $u_1\equiv u_2$. 
 
\rm
 \vskip.15in

 In \eqref{1.3}-\eqref{1.7} we used the notation $f\in H^1(e^{a|x|^{2}}dx)$ if $f,\,\partial_{x_j} f\in L^2(e^{a|x|^{2}}dx)$ for
 $j=1,\dots ,n$, i.e. 
 $$
 \int_{\mathbb  R^n}\,|f(x)|^2\,e^{a|x|^{2}}dx + \,\sum_{j=1}^n \,\int_{\mathbb  R^n}\,|\partial_{x_j} f(x)|^2\,e^{a|x|^{2}}dx< \infty.
$$
 
 By fixing $u_2\equiv 0$ the above  question relates to the persistence properties of solutions of \eqref{1.1c} and  \eqref{1.2c}, i.e.  if $u(x,0)=u_0(x)\in X$ (function space),  then the solution $u=u(x,t)$ 
of \eqref{1.1c}  (resp. \eqref{1.2c}) satisfies that
$$
u\in C([0,T]:X).
$$
 
 These persistence properties, as part of the standard notion of well-posedness, have been studied in function spaces $X$ 
describing the regularity and decay of their elements. For example, they have been established 
 in classical Sobolev spaces $X=H^s(\mathbb R^n), s\geq s_0$, with $s_0$ the optimal Sobolev exponent, which depends on $V$ (resp. $F$), $n$, and $T$, with $T<\infty$, corresponding to local solutions and $T=\infty$ corresponding to global ones, and  in their weighted versions, $X=H^s(\mathbb R^n)\cap L^2(|x|^kdx)$, with $V\in C^{[s+1]}$ (resp. $F\in C^{[s+1]}$),  where persistence holds 
if $s\geq s_0$ and $s \geq k$, due to the fact that $\,\Gamma_j=x_j-2it\partial_{x_j}$, $\,j=1,..,n$, commutes with $\partial_t-i\Delta$, (for details see \cite{HaNaTs} and references therein). In the case, $X=H^s(\mathbb R^n)\cap L^2(|x|^kdx)$, when  $k>s$, persistence fails even in the free case, i.e. $V\equiv 0$ in \eqref{1.1c}, 
and the extra-decay \lq\lq$k-s"$
is transformed  into \lq\lq local regularity".

Note that  \eqref{dada} shows that for the free Schr\"odinger equation, if  $u_0$ is a compactly supported and continuous function, then for any $ t\in \mathbb R\setminus\{0\}$ and $\epsilon >0$, $u(\,\cdot\,, t)$ is not in $L^1(e^{\epsilon|x|}dx)$. In this case, for $t\ne 0$, $u(x, t)$ has an analytic extension to $\mathbb C^n$, so roughly speaking, one can say that the decay, which does not persist  with the solution, is transformed into \lq\lq local regularity".

 One of the key results in  \cite{EsKePoVe}  established that solutions of Schr\"odinger equations with variable coefficients having $L^2$-Gaussian decay at two different times $t_1,\,t_2$, with $t_1<t_2$, preserve this property
in the time interval $[t_1,t_2]$  with a fixed Gaussian weight (see Corollary 2.2 in  \cite{EsKePoVe} and comments after it).
  In this paper, we examine in detail the case of the free particle. We shall  try to understand the possible persistence properties of the solutions of the free Schr\"odinger equation in function spaces with exponential decay at infinity.  

We can summarize part of the results in this work in the following  {\it qualitative} terms: 
\begin{theorem}
\label{Theorem 1.1}   For  $u_0 \in\mathbb S'(\mathbb R^n)$ the following seven statements are equivalent:
\begin{enumerate}
\item[$(i)$] There are two different real numbers $ t_1$ and $t_2$, such that $e^{i t_j\Delta}u_0\in L^2(e^{a_j|x|^2}dx)$, for some  $a_j>0$, $j=1,2$.
\item[$(ii)$] $u_0\in L^2(e^{b_1|x|^2}dx)$ and $\widehat u_0 \in L^2(e^{b_2|x|^2}dx)$, for some $b_j>0$, $j=1,2$.
\item[$(iii)$] There is $\nu: [0,+\infty)\longrightarrow (0,+\infty)$, such that $e^{i t \Delta}u_0\in L^2(e^{\nu(t)|x|^2}dx)$, for all $t\ge 0$.
\item[$(iv)$] $g(x)\equiv e^{i\tau |x|^2}u_0(x)$,  $\tau\in \mathbb R$, verifies $(ii)$ with possible different constants $b_1,\,b_2>0$.

\item[$(v)$] $u_0(x+iy)$ is an entire function such that $\,|u_0(x+iy)|\leq Ne^{-a|x|^2+b|y|^2}$ for some constants $N,\,a,\,b>0$.

\item[$(vi)$] $\widehat u_0(\xi+i\eta)$ verifies $(v)$ with possible different constants $N,\,a,\,b>0$.

\item[$(vii)$] there exist $\delta,\,\epsilon>0$ and $h\in L^2(e^{\epsilon |x|^2}dx)$ such that
$\,u_0(x)=e^{\delta \Delta}h$.
\end{enumerate}
\end{theorem}

Our proof of $(ii)$ implies $(iii)$ in Theorem \ref{Theorem 1.1} will be  a consequence of the new {\it quantitative} results in this work: some  logarithmically convex inequalities for exponentially  weighted  $L^2$-norms of solutions of the free  Schr\"odinger equation. In particular, it will be a  consequence of the following Theorem:

\begin{theorem}
 \label{Theorem 1.2}  Let $\alpha$ and $\beta$ be two positive numbers and  $u$ be the solution of the initial value problem
   \begin{equation*}
   \begin{cases}
   \begin{aligned}
 \label{2.19}
 &\p_t u=i\Delta u,\;\;t>0,\;\;x\in\mathbb R^n,\\
 &u(x,0)=u_0(x).
\end{aligned}
\end{cases}
 \end{equation*}
 Assume that either of the right hand sides of \eqref{2.20}, $\eqref{2.21}$, \eqref{2.22}  or $\eqref{2.23}$ are finite. Then, for any $\lambda\in\mathbb R^n$ and $T>0$ 
 \begin{equation}
\label{2.20}
\|e^{\tfrac{\lambda\cdot x}{\alpha t+\beta}} \,u(t)\|_2
\leq \|e^{\tfrac{\lambda\cdot x}{\beta}} \,u(0)\|_2^{\theta(t)}
\,\|e^{\tfrac{\lambda\cdot x}{\alpha T+\beta}} \,u(T)\|_2^{1-\theta(t)},\;\;\;\;0\leq t\leq T,
\end{equation}
\begin{equation} 
\label{2.21}
\|e^{\tfrac{\lambda\cdot x}{\alpha t+\beta}} \,u(t)\|_2
\leq \|e^{\tfrac{\lambda\cdot x}{\beta}} \,u(0)\|_2^{\mu(t)}
\,\|e^{\frac{2\lambda\cdot \xi}{\alpha}}\, \widehat u(0)\|_2^{1-\mu(t)},
\end{equation}
\begin{equation} 
\begin{aligned}
\label{2.22}
&\;\;\;\;\;\;\,\|e^{\tfrac{|x|^2}{(\alpha t+\beta)^2}} \, u(t)\|_2\\
&\;\;\;\;\;\;\;\;\;\;\;\;\;\;\leq \|e^{\tfrac{|x|^2}{\beta^2}} \,u(0)\|_2^{\theta(t)}\,
\|e^{\tfrac{|x|^2}{(\alpha T+\beta)^2}} \,u(T)\|_2^{1-\theta(t)},\;\;\;0\leq t\leq T,
\end{aligned}
\end{equation}
and
\begin{equation}
\label{2.23}
\|e^{\tfrac{|x|^2}{(\alpha t+\beta)^2}} \,u(t)\|_2
\leq \|e^{\tfrac{|x|^2}{\beta^2}} \,u(0)\|_2^{\mu(t)}\,
\|e^{\frac{4|\xi|^2}{\alpha^2}} \,\widehat u(0)\|_2^{1-\mu(t)},
\end{equation}
where
$$
\theta(t)=\tfrac{\beta(T-t)}{T(\alpha t+\beta)},\;\;\;\;\text{and}\;\;\;\;\mu(t)=\tfrac{\beta}{\alpha t+\beta}.
$$

\end{theorem}

 From (A$'$), (B$'$) and (C$'$), one has that 
 the terms in \eqref{2.23} vanish if $\alpha\beta$ is small enough. The same applies to \eqref{2.22}
for small values of $\alpha$ and $\beta$ depending on $T$. In particular, from (B$'$) in the $1-d$ case, one has that if $4T\ge\beta(\alpha T+\beta)$, then the terms in \eqref{2.22} are zero, and  if $4\ge \alpha \beta$, then the terms in \eqref{2.23} 
 are zero.

In Lemma \ref{Lemma 2.3} we study the class of functions with Gaussian decay and whose Fourier transform also has Gaussian decay. We show that this class  is an algebra with respect to the pointwise product of functions and that it  is closed under the action of the Schr\"odinger group. Therefore, it is also closed with respect to convolutions and multiplication by functions of the form $e^{ia|x|^2},\;a\in\mathbb R$. As a consequence, in Corollary \ref{Corollary 2.1}, we show  that if a free Schr\"odinger solution has Gaussian decay at two different times, then the data belongs to this class, i.e.  
$(i)$  implies $(ii)$ in Theorem \ref{Theorem 1.1}.

Section \ref{S:2} contains the proofs of Theorems \ref{Theorem 1.1} and Theorem \ref{Theorem 1.2}. The proof of Theorem \ref{Theorem 1.2} will be  based on a general abstract result given in Theorem \ref{Theorem 2.1}.

Once Theorem \ref{Theorem 1.2} has been proved we obtain some generalizations and consequences of it. 
Corollaries \ref{Corollary 3.1}-\ref{Corollary 3.3} are extensions of the estimates \eqref{2.22}-\eqref{2.23}. In Corollaries \ref{Corollary 3.4} and \ref{Corollary Viriel interactive identity} we apply  \eqref{2.22}-\eqref{2.23} to any pair (or any finite set) of solutions of the free Schr\"odinger equation.  These estimates, which describe the interaction of two solutions, are somehow similar in spirit to those found in \cite[Sections 10 and 11]{CoKeStTaTa}. Corollary \ref{Corollary 3.5} is an application of the results in Theorem 1.2 and of the Galilean invariance property of the free Schr\"odinger equation. These extensions and applications are in section \ref{S:3}.
  
Finally, we explain and outline some extensions of the results in sections \ref{S:2} and \ref{S:3} to the variable coefficients case. These results will be studied in more detail in a forthcoming publication.

\section{Proof of Theorems \ref{Theorem 1.1} and \ref{Theorem 1.2}}\label{S:2}

We begin with a general abstract result. It can be used to derive properties of logarithmic convexity of  certain $L^2$-norms of solutions of different evolutions.
\begin{theorem}
\label{Theorem 2.1}  Let $\mathcal S$ be a symmetric operator, $\mathcal A$ be an anti-symmetric one, both allowed to depend on the time variable, and $f(x,t)$ be a suitable function. If
 $$
 H(t)=\langle f,f\rangle=\|f(t)\|_2^2,\;\;D(t)=\langle \mathcal S f,f\rangle,\;\;\p_t \mathcal S=\mathcal S_t,\;\;\text{and}\;\;
 N(t)=\frac{D(t)}{H(t)},
 $$
 then
\begin{equation}
 \begin{aligned}
 \label{2.2a}
 &\dot  N(t)=\langle \mathcal S_tf+[\mathcal S,\mathcal A]f,f\rangle /H\\
 &+\frac{1}{2}\Big[ \|\p_tf-\mathcal Af+\mathcal Sf\|_2^2\|f\|_2^2-(\Re\langle \p_tf-\mathcal Af+\mathcal Sf,f\rangle)^2\Big]/H^2\\
 &+\frac{1}{2}\Big[(\Re\langle \p_tf-Af-\mathcal Sf,f\rangle)^2-\|\p_tf-\mathcal Af-\mathcal Sf\|_2^2\|f\|_2^2\Big]/H^2,
 \end{aligned}
 \end{equation}
 and
 \begin{equation*}
 \label{2.3a}
 \dot N(t) \geq\langle \mathcal S_tf+[\mathcal S;\mathcal A]f,f\rangle /H\ -\frac{1}{2}\|\p_tf-\mathcal Af-\mathcal Sf\|_2^2\|f\|_2^2/H^2.
 \end{equation*}
 \end{theorem}
\begin{proof} We have the following identities 
 \begin{equation}
 \label{2.5a}
 \dot H(t)=2\Re\langle \p_tf,f\rangle=\Re\langle\p_tf+\mathcal Sf,f\rangle+\Re\langle\p_tf-\mathcal Sf,f\rangle,
 \end{equation}
  \begin{equation}
 \label{2.6a}
 D(t)=\Re\langle \mathcal Sf,f\rangle=\frac{1}{2}\Re\langle\p_tf+\mathcal Sf,f\rangle-\frac{1}{2}\Re\langle\p_tf-\mathcal Sf,f\rangle.
 \end{equation}
 
 Thus, multiplying \eqref{2.5a} and \eqref{2.6a} it follows that
 \begin{equation}
 \label{2.7a}
 \dot H(t) D(t)=\frac{1}{2}(\Re\langle\p_tf+\mathcal Sf,f\rangle)^2-\frac{1}{2}(\Re\langle\p_tf-\mathcal Sf,f\rangle)^2.
 \end{equation}
 
 Adding an antisymmetric operator does not change the real part, and so
\begin{equation*}
 \label{2.8a}
 \dot H(t) D(t)=\frac{1}{2}(\Re\langle\p_tf-\mathcal Af+\mathcal Sf,f\rangle)^2-\frac{1}{2}(\Re\langle\p_tf-\mathcal Af-\mathcal Sf,f\rangle)^2.
 \end{equation*}
 
 Differentiating $D(t)$
  \begin{equation*}
 \begin{aligned}
 \label{2.9a}
 \dot D(t)&= \langle \mathcal S_tf,f\rangle + \langle \mathcal S\p_tf,f\rangle+\langle \mathcal Sf,\p_tf\rangle\\
 &=\langle \mathcal S_tf,f\rangle + 2\Re \langle \p_tf,\mathcal Sf\rangle\\
 &=\langle \mathcal S_tf+[\mathcal S,\mathcal A]f,f\rangle + 2\Re \langle \p_tf-\mathcal Af,\mathcal Sf\rangle,
 \end{aligned}
\end{equation*}
combined with the polarization identity give
\begin{equation}
 \label{2.10a}
 \dot D(t)=\langle \mathcal S_tf+[\mathcal S,\mathcal A]f,f\rangle +\frac{1}{2}\|\p_tf-\mathcal Af+\mathcal Sf\|_2^2
 -\frac{1}{2}\|\p_tf-\mathcal Af-\mathcal Sf\|_2^2.
 \end{equation}

 The identity \eqref{2.2a} follows from \eqref{2.7a} and \eqref{2.10a}.

Note for future use that 
\begin{equation*}
 \label{2.4a}
 \dot H(t)=2\Re \langle \p_tf-\mathcal Sf-\mathcal Af,f\rangle +2D(t),
 \end{equation*}
 so if  $\p_tf-\mathcal Sf-\mathcal Af=0$, one also has that $\dot H(t)=2D(t)$.
  \end{proof}
 
\begin{remark} \label{megusta}

The abstract identity in Theorem \ref{Theorem 2.1} shows in disguise  the {\it \lq\lq frequency function\rq\rq} or \lq\lq{\it monotonicity argument}\rq\rq\  linked to the Carleman inequality
\[\|\partial_tf-\mathcal Af\|^2_2+ \|\mathcal Sf\|^2_2\le \|\partial_tf-\mathcal Sf-\mathcal Af\|^2_2\ .\]

The antisymmetric and symmetric parts of $\partial_t-\mathcal S-\mathcal A$, as a space-time operator, are respectively $\partial_t-\mathcal A$ and $-\mathcal S$. Its commutator, $[-\mathcal S,\partial_t-\mathcal A]$, is $\mathcal S_t+[\mathcal S,\mathcal A]$. Thus,
\begin{align*}
 \|\partial_tf-\mathcal Sf-\mathcal Af\|^2_2&=\|\partial_tf-\mathcal Af\|^2_2+ \|\mathcal Sf\|^2-2\text{\it Re} \iint \mathcal Sf\overline{\partial_t f-\mathcal Af}\,dxdt\\
 =&\|\partial_tf-\mathcal Af\|^2_2+ \|\mathcal Sf\|^2_2+\iint[-\mathcal S,\partial_t-\mathcal A]f\overline f\,dxdt\\
 =&\|\partial_tf-\mathcal Af\|^2_2+ \|\mathcal Sf\|^2_2+\iint \left(\mathcal S_tf+[\mathcal S, \mathcal A]f\right)\overline f\,dxdt\ ,
\end{align*}
and the Carleman inequality holds, when $\mathcal S_t+[\mathcal S, \mathcal A]$ is non-negative. Theorem \ref{Theorem 2.1} shows that $H(t)$ is logarithmically convex, when $\mathcal S_t+[\mathcal S, \mathcal A]$ is non-negative, $\partial_tf-\mathcal Sf-\mathcal Af=0$, and provided that the calculations and integrations by parts carried out in the application of Theorem \ref{Theorem 2.1} to a particular case can be justified. 
\end{remark}
 We apply Theorem \ref{Theorem 2.1} in our proof of Theorem \ref{Theorem 1.2}, and in order to justify the finiteness of the quantities involved, integrations by parts or calculations involved in this application of Theorem \ref{Theorem 2.1}, we use Lemma \ref{L: StSh} \cite[pp.130]{StSh}:
 
\begin{lemma}\label{L: StSh} Let $f$ be an entire function such that
\begin{equation}\label{E: clasedefunciones}
|f(x+iy)|\le Ne^{-a|x|^2+b|y |^2},\;\,\text{with}\;\,N,\,a, b>0,\;\;\forall \,x,\,y\in\mathbb R^n.
\end{equation}
Then, $\widehat f$ is an entire function and  
\begin{equation*}
|\widehat f(\xi+i\eta)|\le N'e^{-a'|\xi|^2+b'|\eta|^2},\;\,\;\;\forall \,\xi,\,\eta\in\mathbb R^n\ 
\end{equation*}
for some positive constants $N'$, $a'$ and $b'$.
\end{lemma}

\begin{proof}[Proof of Theorem \ref{Theorem 1.2}]  We apply Theorem \ref{Theorem 2.1} with
\begin{equation}\label{E: estafuncion}
f(x,t)=e^{\frac{\lambda\cdot x}{\alpha t+\beta}} u(x,t),\;\;\;\;\lambda\in\mathbb R^n.
\end{equation}

When the initial data to the free Schr\"odinger equation verifies that the right hand side of  \eqref{2.23} is finite, $u_0$ is in  $H^\infty(\mathbb R^n)$ and $u$ is in $C^\infty(\mathbb R:H^\infty(\mathbb R^n))$. Also, $u_0$ extends to the complex-space $\mathbb C^n$ as an analytic function, and there are positive constants $N$, $a$ and $b$, such that
\begin{equation*}
|u_0(\xi+i\eta)|\le Ne^{-a|\xi|^2+b|\eta |^2}\ \;\;\text{for all}\ \xi, \eta\in\mathbb R^n\ .
\end{equation*}

When $T$ is positive, $f(z)=e^{\frac{iz^2}{4T}}u_0(z)$, verifies the conditions in Lemma \ref{L: StSh} and $u(T)$ is essentially the Fourier transform of $e^{\frac{i|y|^2}{4T}}u_0(y)$. Thus, 
$$
\|e^{a'|x|^2}u(T)\|_2
$$
 is finite for some positive number $a'$. Corollary 2.2 in \cite{EsKePoVe} shows that,
 $$
 \sup_{0\le t\le T}\|e^{a''|x|^2}u(t)\|_2<+\infty ,
 $$
for some new $a''$, which might depend on $a$ and $b$.  The formula \eqref{dada} shows that the same occurs, when  the right hand side of \eqref{2.22} is finite. Moreover, the same holds, in both cases, for all the derivatives of $u$.

Once this has been settled, our choice of $f$ in \eqref{E: estafuncion} shows that
$$
\p_t f = \mathcal Sf +\mathcal Af,
$$
where
$$
\mathcal S=-\frac{2 i}{\alpha t+\beta}\,\lambda\cdot \nabla -\frac{\alpha}{(\alpha t+\beta)^2}\,\lambda\cdot x,\;\;\;\;
\mathcal A=i(\Delta +\frac{|\lambda|^2}{(\alpha t+\beta)^2}),
$$
$$
\mathcal S_t+[\mathcal S,\mathcal A]=- \frac{2 \alpha}{\alpha t+\beta} \mathcal S
$$
and if
 $$
 H(t)=\int_{\mathbb R^n}|f(x,t)|^2\,dx=\int_{\mathbb R^n} e^{\frac{2\lambda\cdot x}{\alpha t+\beta}} |u(x,t)|^2\,dx,
 $$

 The last comment in the poof of Theorem \ref{Theorem 2.1} shows that
 $$
 \p_t \log H(t)=2 N(t) =2\,\frac{\langle Sf,f\rangle}{H}
 $$
 and
 \begin{equation}
 \begin{aligned}
 \label{2.24}
  \p^2_t \log H(t)&=2 \dot N(t) \geq 2\langle S_tf+[S;A]f,f\rangle/H\\
  &\geq-\frac{4\alpha}{\alpha t+\beta}\frac{\langle Sf,f\rangle}{H}
  = - \frac{2\alpha}{\alpha t+\beta} \p_t \log H(t).
 \end{aligned}
 \end{equation}
 In particular, \eqref{2.24} implies that the function
  $$
  G(t)=H(t)^{\alpha t+\beta},\;\;\;\;\;0\leq t\leq T,
  $$
  is logarithmically  convex. Thus,
  $$
  G(t)\leq G(0)^{(T-t)/T}\,G(T)^{t/T},\;\;\;\;\;\;0\leq t\leq T,
  $$
  and consequently
 \begin{equation}\label{E: linealcontiempo}
 H(t)\leq H(0)^{\beta (T-t)/T(\alpha t+\beta)}\,H(T)^{t(\alpha T+\beta)/T(\alpha t+\beta)},
\end{equation}
  which yields \eqref{2.20}.
  
  To prove \eqref{2.21} we recall the formulae
  \begin{equation}
  \begin{aligned}
  \label{a.1}
  u(x,t)&=(2\pi)^{-n/2}\int_{\mathbb R^n}e^{i x\cdot \xi-it|\xi|^2}\widehat u_0(\xi)d\xi
  \\
  &=(4\pi i t)^{-n/2}\int_{\mathbb R^n}e^{\frac{i|x-y|^2}{4t}}u_0(y)dy.
 \end{aligned}
  \end{equation}
  Thus
 \begin{equation*}
 \begin{aligned}
  \label{a.2} 
 u\left(\frac{x}{t},\frac{1}{t}\right)=&(2\pi)^{-n/2}\int_{\mathbb R^n}e^{\frac{i x\cdot \xi}{t}-\frac{i|\xi|^2}{t}}\widehat u_0(\xi)d\xi\\
  =&(2\pi )^{-n/2}\int_{\mathbb R^n}e^{\frac{i|x-2\xi|^2}{4t}+\frac{i|x|^2}{4t}}\widehat u_0(\xi)d \xi,
  \end{aligned}
 \end{equation*}
  and so
   \begin{equation}
   \begin{aligned}
  \label{a.3}
  (-it)^{-n/2}&e^{-\frac{i|x|^2}{4t}}u\left(\frac{x}{t},\frac{1}{t}\right)\\
  &\;\;\;\;\;\;\;\;\;= (4\pi i t)^{-n/2}\int_{\mathbb R^n}e^{\frac{-i|x-\xi|^2}{4t}}2^{-n/2}\, \widehat u_0(\xi/2)d \xi.
\end{aligned}
  \end{equation}
  
  Hence, using the $\psi$-conformal or Appel transformation we define
 \begin{equation}
 \label{2.29}
\overline v(x,t)=(-it)^{-n/2}e^{-\frac{i|x|^2}{4t}}u\left(\frac{x}{t},\frac{1}{t}\right),
\end{equation}
and see from \eqref{a.1}-\eqref{a.3} that $v(x,t)$ is the solution of the initial value problem
 \begin{equation}
   \begin{cases}
   \begin{aligned}
 \label{2.30}
 &\p_t v=i\Delta v,\;\;t>0,\;\;x\in\mathbb R^n,\\
 &v(x,0)=2^{-n/2}\,\overline{\widehat u_0}(x/2).
\end{aligned}
\end{cases}
 \end{equation}

Assume now that $0<t<1$. From \eqref{2.20}, with $T=1$, it follows that
\begin{equation}
\label{2.31}
\|e^{\frac{\lambda\cdot x}{\alpha t+\beta}} u(t)\|_2
\leq \|e^{\frac{\lambda\cdot x}{\beta}}u(0)\|_2^{\frac{\beta(1-t)}{\alpha t+\beta}}
\|e^{\frac{\lambda\cdot x}{\alpha +\beta}}u(1)\|_2^{\frac{(\alpha +\beta)t}{\alpha t+\beta}}.
\end{equation}

 Interchanging the role of $\alpha$ and $\beta$ in  \eqref{2.20} and applying it with 
 $T=1/t>1$ to $v(x,t)$, the solution of \eqref{2.30}, one gets 
\begin{equation}
\label{2.32}
\|e^{\frac{\lambda\cdot x}{\beta +\alpha}} v(1)\|_2
\leq \|e^{\frac{\lambda\cdot x}{\alpha}}v(0)\|_2^{\frac{\alpha(1-t)}{\alpha +\beta}}
\|e^{\frac{\lambda\cdot x}{\alpha +\beta t}}v(1/t)\|_2^{\frac{(\beta+\alpha t)}{\alpha +\beta}}.
\end{equation}

Since from \eqref{2.29}
\begin{equation*}
\label{2.33}
\|e^{\frac{\lambda\cdot x}{\beta +\alpha}} v(1)\|_2=\|e^{\frac{\lambda\cdot x}{\beta +\alpha}} u(1)\|_2,
\end{equation*}
and
\begin{equation}
\label{2.34}
\|e^{\frac{\lambda\cdot x}{\beta/t +\alpha }} v(1/t)\|_2=\|e^{\frac{\lambda\cdot x}{\alpha t+\beta}} u(t)\|_2,
\end{equation}
combining \eqref{2.32}-\eqref{2.34} and the value of $v$ at the initial time, we get
\begin{equation}
\label{2.35}
\|e^{\frac{\lambda\cdot x}{\alpha +\beta}} u(1)\|_2
\leq \|e^{\frac{2\lambda\cdot \xi}{\alpha}}\widehat u(0)\|_2^{\frac{\alpha(1-t)}{\alpha +\beta}}
\|e^{\frac{\lambda\cdot x}{\alpha t +\beta}}u(t)\|_2^{\frac{\alpha t+\beta}{\alpha +\beta}}.
\end{equation}

From \eqref{2.31}, and \eqref{2.35}, and the fact that $\,\|e^{\frac{\lambda\cdot x}{\alpha t +\beta}}u(t)\|_2$ is finite, we get 
 \begin{equation}
 \label{arriba}
 \|e^{\frac{\lambda\cdot x}{\alpha t+\beta}} u(t)\|_2
\leq \|e^{\frac{\lambda\cdot x}{\beta}}u(0)\|_2^{\frac{\beta}{\alpha t+\beta}}
\|e^{\frac{2 \lambda\cdot \xi}{\alpha}}\widehat u(0)\|_2^{\frac{ \alpha t}{\alpha t+\beta}},
\end{equation}
which gives  \eqref{2.21}.

To prove \eqref{2.23} we square both sides of \eqref{arriba}, multiply them by
$\,e^{-|\lambda|^2/2}$, integrate with respect to $ \,\lambda\in\mathbb R^n$, and use the identity
\begin{equation}
\label{key1}
\int_{\mathbb R^n}e^{\frac{2\lambda\cdot x}{\gamma}-\frac{|\lambda|^2}{2}}\,d\lambda =(2\pi)^{n/2}
e^{\frac{2|x|^2}{\gamma^2}},\;\;\;\;\;\;\gamma\in\mathbb R,
\end{equation}
to obtain in the left hand side
\begin{equation}
\begin{aligned}
\label{explain1}
&\int_{\mathbb R^n} e^{-\frac{|\lambda|^2}{2}}\|e^{\frac{\lambda\cdot x}{\alpha t+\beta}} u(t)\|_2^2d\lambda =
\int_{\mathbb R^n}\int_{\mathbb R^n} e^{-\frac{|\lambda|^2}{2}} e^{\frac{2\lambda\cdot x}{\alpha t+\beta}} |u(x,t)|^2dxd\lambda\\
&= \int_{\mathbb R^n} (\int_{\mathbb R^n} e^{\frac{2\lambda\cdot x}{\alpha t+\beta}-\frac{|\lambda|^2}{2}}\,d\lambda) |u(x,t)|^2dx\\
& = (2\pi)^{n/2} 
\int_{\mathbb R^n} e^{\frac{2 |x|^2}{(\alpha t+\beta)^2}} |u(x,t)|^2dx,
\end{aligned}
\end{equation}
and on the right hand side a combination of  a similar argument to that in \eqref{explain1} with   H\"older inequality,
$p=\frac{\alpha t+\beta}{\beta},\,p'=\frac{\alpha t+\beta}{\alpha t}$ leads to
\begin{equation}
\begin{aligned}
\label{explain2}
&\int_{\mathbb R^n} e^{-\frac{|\lambda|^2}{2}} \Big(\int e^{\frac{2\lambda\cdot x}{\beta}} 
|u(x,0)|^2 dx\Big)^{\frac{1}{p}} \Big(\int e^{\frac{4\lambda\cdot \xi}{\alpha}} 
|\widehat u(\xi,0)|^2 d\xi\Big)^{\frac{1}{p'}}\;d\lambda\\
&\le\Big(\int \int e^{-\frac{|\lambda|^2}{2}+\frac{2\lambda\cdot x}{\beta}} |u(x,0)|^2dx d\lambda\Big)^{\frac{\beta}{\alpha t+\beta}}\\
&\;\;\;\cdot\; 
\Big(\int \int e^{-\frac{|\lambda|^2}{2}+\frac{4\lambda\cdot \xi}{\alpha}} |\widehat u(\xi,0)|^2d\xi d\lambda\Big)^{\frac{\beta}{\alpha t+\beta}}\\
&=(2\pi)^{n/2} \Big(\int_{\mathbb R^n} e^{\frac{2 |x|^2}{\beta^2}} |u(x,0)|^2dx\Big)^{\frac{2\beta}{\alpha t+\beta}}\\
&\,\;\;\;\;\;\cdot \Big(\int_{\mathbb R^n} e^{\frac{8 |\xi|^2}{\alpha^2}} |\widehat u(\xi,0)|^2d\xi\Big)^{\frac{\alpha t}{\alpha t+\beta}}.
\end{aligned}
\end{equation}
 Therefore, \eqref{explain1} and \eqref{explain2} yield \eqref{2.23}.

Finally, to obtain \eqref{2.22} we reapply the  argument  used to get \eqref{explain1}-\eqref{explain2} but starting with \eqref{2.20} instead of \eqref{2.21}. 
\end{proof}

Next, we introduce some notation: 
for $f\in\mathbb S'(\mathbb R^n)$, $\,p,\,q\in[1,\infty]$, and $A_1, A_2>0$ we will write  
\begin{equation*}
\label{2.50}
f=O_{p,q}(A_1;A_2)\;\text{if}\; e^{A_1|\cdot|^2}f\in L^p(\mathbb R^n),\,e^{A_2|\cdot|^2}\widehat f\in L^q(\mathbb R^n),\\
\;\text{with}\;M_{p,q}(f)=A_1A_2.
\end{equation*}

In the case $p=q$ we shall write $\;O_{p}$ instead of $O_{p,p}$.

In this context, Hardy's result mentioned in the introduction and its extension to higher dimension found in \cite{SiSu} tells us that  if $M_{\infty,\infty}(f)> 1/4$, then $f\equiv 0$. Also a result in \cite{CoPr}
affirms that in the 1-d case if $M_{p,q}(f)\geq  1/4$, for any $p,\,q\in[1,\infty]$ with at least one of them finite, then $f\equiv 0$.

\begin{lemma} \label{Lemma 2.3} 
  If $f, \,g,\, h \in\mathbb S'(\mathbb R^n)$ are such that $f=O_{p,q}(A_1;A_2),\;g=O_{r,q}(B_1;B_2)$
 and $h=O_{2}(C_1;C_2)$, with $p, q, r \in[1,\infty]$ and $1/p+1/r=1/l$, then
\begin{itemize}
\item[(a)] $\widehat f=O_{q,p}(A_2;A_1)$.
\item[(b)] $fg = O_{l,q}(A_1+B_1;A_2B_2/(A_2+B_2))$ and
 \[M_{l,q}(fg)\leq \frac{B_2}{A_2+B_2}M_{p,q}(f)+ \frac{A_2}{A_2+B_2} M_{r,q}(g).\]
 \item[(c)] $g\ast f = O_{q,l}(A_1B_1/(A_1+B_1); A_2+B_2)$.
 \item[(d)] $e^{it\Delta} h=O_2(C_1C_2/(C_2+4t\sqrt{C_1C_2}+4t^2C_1);\,C_2)$, when $t>0$.
 \item[(e)] $e^{-i\tau |\cdot|^2} h = O_2(C_2; \,C_1C_2/(C_1+4 \tau\sqrt{C_1C_2}+4\tau^2C_2))$, when $\tau>0$.
\end{itemize}  
\end{lemma}
\begin{proof} Part (a) is immediate. To obtain (b) we need the following calculation:
for $\mu, \nu>0$
\begin{equation}
\label{calculo}
\begin{aligned}
\int_{\mathbb R^n}e^{-\mu|x-y|^2}e^{-\nu|y|^2}dy
&= e^{-(\mu-\frac{\mu^2}{\mu+\nu})|x|^2}
\int_{\mathbb R^n}e^{-|\frac{\mu x}{\sqrt{\mu+\nu}}-\sqrt{\mu+\nu}\, y|^2}dy\\
&=(\tfrac{\pi}{\mu+\nu})^{n/2}e^{-\frac{\mu\nu}{\mu+\nu}|x|^2}.
\end{aligned}
\end{equation}

Then, taking $c=A_2B_2/(A_2+B_2)$, and combining  \eqref{calculo} and H\"older's inequality we get
$$
\begin{aligned}
&\|e^{c|x|^2}\widehat f\ast \widehat g\|_q
=\|e^{c|x|^2}\int e^{A_2|x-y|^2}\widehat f(x-y) e^{B_2|y|^2} \widehat g(y) e^{-A_2|x-y|^2}e^{-B_2|y|^2}dy\|_q\\
&\leq \|e^{c|x|^2}(\int e^{-q'A_2|x-y|^2}e^{-q'B_2|y|^2}dy)^{1/q'}\|_{\infty}\\&\;\;\;\;\;\;\;\;\;\;
\|(\int e^{qA_2|x-y|^2}|\widehat f(x-y)|^q e^{qB_2|y|^2} |\widehat g(y)|^q dy)^{1/q}\|_q\\
&\leq c \| e^{A_2|x|^2}\,\widehat f\,\|_q\,\|e^{B_2|x|^2} \,\widehat g\|_q,
\end{aligned}
$$
which proves (b).

A combination of (a) and (b) yields (c). 
Part (d) follows from Theorem \ref{Theorem 1.2} estimate \eqref{2.23}. Finally, (e) follows by combining (a), (d) and the formula
\[
e^{-i\tau |x|^2} h(x)=\widehat{ e^{i\tau \Delta} (h^{\lor})}(x).
\]
\end{proof}
\begin{remark}\label{Unanota} 
Parts (d) and (e) in Lemma \ref{Lemma 2.3}  still hold with  $t,\, \tau\in \mathbb R$ 
resp.  by replacing $t$ and $\tau$ on their right hand side by $|t|$ and $|\tau|$ resp.
\end{remark}

\begin{remark}\label{R: remarkdeluis}
The Lemma suggests to consider the class of funtions verifying \eqref{E: clasedefunciones} as a space of test functions for more general distributions than the tempered ones. This issue will be studied elsewhere.
\end{remark}

\begin{corollary}
 \label{Corollary 2.1}
  Let $u_0\in\mathbb S'(\mathbb R^n)$. If 
  \begin{equation*}
  \label{2.51}
  e^{it_j\Delta} u_0 \in L^2(e^{2\mu_j |x|^2}dx),\;j=1,2,\;\;t_1\neq t_2,\;\;\mu_1,\,\mu_2>0,
  \end{equation*}
  then
  \begin{equation*}
  \label{2.52}
  u_0\in O_2(A_1;A_2),\;\;\;\text{for some}\;\;A_1,\,A_2>0.
  \end{equation*}
  \end{corollary}
 
\begin{proof} Using part (d) in Lemma \ref{Lemma 2.3} and Remark \ref{Unanota} we can assume $t_1=0, \,t_2=s>0$. Thus, 
 \begin{equation}
  \label{2.54}
 f=u_0 \in L^2(e^{2\mu_1 |x|^2}dx),\;\;\text{and}\;\;\;e^{i s\Delta}f\in L^2(e^{2\mu_2 |x|^2}dx).
  \end{equation}
  But
  $$
  \aligned
  e^{is\Delta}f(x)&=(4\pi i s)^{-n/2}\int_{\mathbb R^n}e^{\frac{i|x-y|^2}{4s}}f(y)dy\\
  & =(4\pi i s)^{-n/2}
 e^{\frac{i|x|^2}{4s}}\int_{\mathbb R^n} e^{\frac{- i x\cdot y}{2s}}e^{\frac{i|y|^2}{4s}}
 f(y)dy.
 \endaligned
 $$
 
 Therefore,
 \begin{equation}
  \label{2.55}
e^{\frac{-i|x|^2}{4s}} e^{i s\Delta}f(x) =(4\pi i s)^{-n/2}\,\widehat{e^{\frac{i|\cdot|^2}{4s}}
 f}(x/s).
 \end{equation}

From \eqref{2.54}-\eqref{2.55} it follows that
 $$
 e^{\frac{i|\cdot|^2}{4s}} f\in L^2(e^{2\mu_1 |x|^2}dx),\;\;\text{and}\;\;\;
 \widehat{e^{\frac{i|\cdot|^2}{4s}}
 f} \in L^2(e^{2\mu_2(2s)^2 |\xi|^2}d\xi).
 $$
In particular, $e^{\frac{i|\cdot|^2}{4s}} f=O_2(\mu_1; 4s^2 \mu_2)$, and from Lemma \ref{Lemma 2.3}, part (e) and the Remark \ref{Unanota}, it follows that
$$
u_0=f=O_2(4s^2\mu_2 ; \frac{4 s^2\mu_1 \mu_2}{
\mu_1+2\sqrt{\mu_1\mu_2}+\mu_2}),
$$
which yields the desired result.
\end{proof}

 \begin{proof}[Proof of Theorem \ref{Theorem 1.1}] The part $(i)$ implies $(ii)$ is Corollary \ref{Corollary 2.1}. That $(ii)$ implies $(iii)$ is in 
 Theorem \ref{Theorem 1.2}, \eqref{2.23}. $(iii)$ implies $(i)$ is immediate. Lemma \ref{Lemma 2.3} part $(e)$ shows that $(i)$ and $(iv)$ are equivalent. Lemma \ref{L: StSh} affirms that $(v)$ and $(vi)$ are equivalent and that each one implies $(ii)$. The fact that $(i)$ implies $(vi)$ follows from the following result (see \cite{StSh}, page 130):
 
 \it Suppose that $f(x+iy)$ defined in $\mathbb C^n$ is an entire function such that $f(z)=O(e^{c_1|z|^2})$ for some $c_1>0$ with $f(x)=O(e^{-c_2|x|^2})$ for some $c_2>0$. Then 
 $|f(x+iy)|=O(e^{-a|x|^2+b|y|^2})$ for some $a,\,b>0$. \rm

  \vskip.05in
\noindent $(vii)$ implies $(ii)$ is immediate using Lemma \ref{Lemma 2.3} part $(b)$ .  To see that $(ii)$ implies $(vii)$ define 
 $$
 \widehat h(\xi)=e^{\delta|\xi|^2} \widehat u_0(\xi)\in L^1(\mathbb R^n) \cap L^{\infty}(\mathbb R^n),\;\;\;\;\delta\in (0,b_2),
 $$
 So we just need to show that $h\in L^2(e^{\epsilon|x|^2}dx)$ for some $\epsilon>0$. Using that $(ii)$ implies $(vi)$ it follows that
 $$
 |\widehat h(\xi+i\eta)|=|e^{\delta(\xi+i\eta)\cdot(\xi+i\eta)} \widehat u_0(\xi+i\eta)|
 \leq N e^{(-a+\delta)|\xi|^2+b_1|\eta|^2},\;\;\;\;\;a,\,b_1>0.
 $$
 Hence, taking $\delta<a$ and using that $(v)$ and $(vi)$ are equivalent we get that $h\in L^2(e^{\epsilon|x|^2}dx)$ for some $\epsilon>0$ which completes the proof.
  \end{proof}

\section{Further results. Generalizations and Applications}\label{S:3}
\subsection{Some other convex weights}
We return to Theorem \ref{Theorem 1.1}, and its proof given in section \ref{S:2}. From the arguments used in \eqref{2.22}-\eqref{2.23}
and \eqref{explain1}-\eqref{explain2}, it is  clear  that similar estimates hold with different  coefficients multiplying the Gaussian weight in each variable. More precisely, we have the following result:

\begin{corollary}
\label{Corollary 3.1}  Using the same hypotheses and notation that in Theorem \ref{Theorem 1.2}. 
Given $\vec \gamma =(\gamma_1,\dots,\gamma_n)\in [0,\infty)^n$
(using summation convention over multiple indices) one has that
\begin{equation*} 
\label{2.22b}
\|e^{\frac{\gamma_j x_j^2}{(\alpha t+\beta)^2}} u(t)\|_2
\leq \|e^{\frac{\gamma_j x_j^2}{\beta^2}}u(0)\|_2^{\theta(t)}\,
\|e^{\frac{\gamma_j x_j^2}{(\alpha T+\beta)^2}}u(T)\|_2^{1-\theta(t)},
\end{equation*}
for any $0\le t\le T$, and
\begin{equation*}
\label{2.23b}
\|e^{\frac{\gamma_j x_j^2}{(\alpha t+\beta)^2}} u(t)\|_2
\leq \|e^{\frac{\gamma_j x_j^2}{\beta^2}}u(0)\|_2^{\mu(t)}\,
\|e^{\frac{4\gamma_j x_j^2}{\alpha^2}}\widehat u(0)\|_2^{1-\mu(t)},
\end{equation*}
when $t\ge 0$.
\end{corollary}

Next, we shall extend Theorem \ref{Theorem 1.2} to the case where we replace the quadratic  powers in the exponents in 
\eqref{2.22}-\eqref{2.23} with possibly  different powers in each component. 

\begin{corollary}
 \label{Corollary 3.2} 
 Using the same hypotheses and notation that in Theorem 1.2. 
 Given  $\vec p=(p_1,..,p_n)\in (1,2]^n$ and  $\vec \gamma =(\gamma_1,\dots,\gamma_n)\in [0,\infty)^n$ there exists $c=c(\vec p)>0$ such that
\begin{equation} 
\label{lpversion1}
\|e^{\gamma_j\left|\frac{x_j}{\alpha t+\beta}\right|^{p_j}} u(t)\|_2
\leq c \; \|e^{\gamma_j \left|\frac{x_j}{\beta}\right|^{p_j}} u(0)\|_2^{\theta(t)}\;\|e^{\gamma_j \left|\frac{x_j}{\alpha T+\beta}\right|^{p_j}} u(T)\|_2^{1-\theta(t)},
\end{equation}
for $0\le t\le T$, and 
\begin{equation}
\label{lpversion2}
\|e^{\gamma_j \left|\frac{x_j}{\alpha t+\beta}\right|^{p_j}}u(t)\|_2
\leq c \;\|e^{\gamma_j \left|\frac{x_j}{\beta}\right|^{p_j}} u(0)\|_2^{\mu(t)}\;
 \|e^{\gamma_j\,\left|\frac{2\xi_j}{\alpha} \right|^{p_j}} \widehat u(0)\|_2^{1-\mu(t)},
 \end{equation}
 for $t\ge 0$.
 \end{corollary}
\begin{proof} First, we notice that instead of the identity \eqref{key1} one has the following asymptotic formula (see \cite{StSh} Proposition 2, pp. 323):  in the one dimensional case 
 \begin{equation*}
\int_{\mathbb R}e^{\lambda x-\frac{|\lambda|^{p'}}{p'}}|\lambda|^{\frac{p'-2}2}\,d\lambda=e^{\frac{|x|^p}p}\left(\left(2\pi\right)\sqrt{p-1}+\text{O}\left(|x|^{-\frac{p}2}\right)\right),
\end{equation*}
when $1<p<\infty$, $|x|\ge 1$ and $1/p+1/p'=1$. In particular, there is $c=c(p)$ such that 
 \begin{equation}
 \label{form1}
c^{-1}e^{\frac{|x|^p}p}\le \int_{\mathbb R}e^{\lambda x-\frac{|\lambda|^{p'}}{p'}}|\lambda|^{\frac{p'-2}2}\,d\lambda\le c e^{\frac{|x|^p}p},
\end{equation}
when $x$ is in $\mathbb R$. Thus, to obtain \eqref{lpversion1}-\eqref{lpversion2}, one just  follows the  argument provided in the proof of Theorem \ref{Theorem 1.2} to obtain \eqref{2.22} and \eqref{2.23} respectively, but using instead of the identity \eqref{key1}, the inequality \eqref{form1}, when $n=1$ and \eqref{form2}, when $n\ge 2$:
\begin{equation}
\label{form2}
c^{-1}e^{\frac{|x|^p}p}\le \int_{\mathbb R^n}e^{\lambda\cdot x-\frac{|\lambda|^{p'}}{p'}}|\lambda|^{\frac{n(p'-2)}2}\,d\lambda\le c\,e^{\frac{|x|^p}p},
\end{equation}
when $1<p<\infty$ and $|x|\ge 1$.
\end{proof}
\begin{corollary}
 \label{Corollary 3.3} With the same hypotheses and notation as in Theorem \ref{Theorem 1.2}. Given any $p\in (1,2]$, there is $c=c(p)>0$, such that
\begin{equation*} 
\label{lpversion3}
\|e^{\left|\frac{x}{\alpha t+\beta}\right|^{p}} u(t)\|_2
\leq c \; \|e^{\left|\frac{x}{\beta}\right|^{p}} u(0)\|_2^{\theta(t)}\;\|e^{\left|\frac{x}{\alpha T+\beta}\right|^{p}} u(T)\|_2^{1-\theta(t)},
\end{equation*}
for $0\le t\le T$, and 
\begin{equation*}
\label{lpversion4}
\|e^{\left|\frac{x}{\alpha t+\beta}\right|^{p}}u(t)\|_2
\leq c \;\|e^{\left|\frac{x}{\beta}\right|^{p}} u(0)\|_2^{\mu(t)}\;
 \|e^{\left|\frac{2\xi}{\alpha} \right|^{p}} \widehat u(0)\|_2^{1-\mu(t)},
 \end{equation*}
for $t\ge 0$.
\end{corollary}

 Note that we have stated these results for $1<p\leq 2$ since Hardy's uncertainty principle shows that for $p>2$ all the functions are $0$.
 
\subsection{Products of solutions} Next, we shall apply in Corollary \ref{Corollary 3.4} the logarithmically convex inequalities to a pair of solutions $e^{it\Delta}u_0$ and $e^{it\Delta}v_0$: starting with the inequality \eqref{2.20} or \eqref{2.21} for each of these solutions, multiplying the squares of their left hand sides and their right ones respectively, and following the argument in \eqref{key1}, \eqref{explain1} and \eqref{explain2}, we obtain some estimates concerning the interaction of two solutions of the free Schr\"odinger equation. 
\begin{corollary}
 \label{Corollary 3.4} Under the same hypotheses and notation that in Theorem \ref{Theorem 1.2}, the following inequalities hold
 
\begin{multline}\label{dos1}
\|e^{\frac{\lambda\cdot(x-y)}{\alpha t+\beta}} e^{it\Delta}u_0(x)\,e^{it\Delta}v_0(y)\|_{L^2(\mathbb R^{2n}_{x,y})}\\
\leq \|e^{\frac{\lambda\cdot(x-y)}{\beta}}u_0(x)\,v_0(y)\|_{L^2(\mathbb R^{2n}_{x,y})}^{\theta(t)}\,
\|e^{\frac{\lambda\cdot(x-y)}{\alpha T+\beta}} e^{iT\Delta} u_0(x)\, e^{iT\Delta}v_0(y)\|_{L^2(\mathbb R^{2n}_{x,y})}^{1-\theta(t)},
\end{multline}
\begin{multline}\label{dos3}
\|e^{\frac{|x-y|^2}{(\alpha t+\beta)^2}} e^{it\Delta}u_0(x)\,e^{it\Delta}v_0(y)\|_{L^2(\mathbb R^{2n}_{x,y})}\\
\leq \|e^{\frac{|x-y|^2}{\beta^2}}u_0(x)\,v_0(y)\|_{L^2(\mathbb R^{2n}_{x,y})}^{\theta(t)}\,
\|e^{\frac{|x-y|^2}{(\alpha T+\beta)^2}} e^{iT\Delta}u_0(x)\,e^{iT\Delta} v_0(y)\|_{L^2(\mathbb R^{2n}_{x,y})}^{1-\theta(t)},
\end{multline}
for $0\le t\le T$, and
\begin{multline*}
\|e^{\frac{\lambda\cdot(x-y)}{\alpha t+\beta}} e^{it\Delta}u_0(x)\,e^{it\Delta}v_0(y)\|_{L^2(\mathbb R^{2n}_{x,y})}\\
\leq \|e^{\frac{\lambda\cdot(x-y)}{\beta}}u_0(x)\,v_0(y)\|_{L^2(\mathbb R^{2n}_{x,y})}^{\mu(t)}\,
\|e^{\frac{2\lambda\cdot(\xi-\eta)}{\alpha}}\widehat u_0(\xi)\,\widehat v_0(\eta)\|_{L^2(\mathbb R^{2n}_{\xi,\eta})}^{1-\mu(t)},
\end{multline*}
\begin{multline*} 
\|e^{\frac{|x-y|^2}{(\alpha t+\beta)^2}} e^{it\Delta}u_0(x)\,e^{it\Delta}v_0(y)\|_{L^2(\mathbb R^{2n}_{x,y})}\\
\leq \|e^{\frac{|x-y|^2}{\beta^2}}u_0(x)\,v_0(y)\|_{L^2(\mathbb R^{2n}_{x,y})}^{\mu(t)}
\|e^{\frac{4|\xi-\eta|^2}{\alpha^2}}\widehat u_0(\xi)\,\widehat v_0(\eta)\|_{L^2(\mathbb R^{2n}_{\xi,\eta})}^{1-\mu(t)},
\end{multline*}
for $t\ge 0$.
\end{corollary}

\begin{remark}\label{Notados}
The interaction inequalities \eqref {dos1}-\eqref{dos3} show that the $\alpha t+\beta$ power of their left hand sides are logarithmically convex functions in $[0,T]$. 

  Also observe that appropriate versions of these inequalities can be deduced for any finite set of solutions of the free Schr\"odinger equation.
\end{remark}
In Corollary \ref{Corollary Viriel interactive identity}, we show that the \lq\lq logarithmic convexity\rq\rq\  behind the interaction inequalities in Corollary \ref{Corollary 3.4}, implies interaction Morawetz inequalities for the free particles in the same  spirit as the interaction inequalities in \cite[Sections 10 and 11]{CoKeStTaTa}.
\begin{corollary}\label{Corollary Viriel interactive identity}  Under the same hypotheses and notation than in  Corollary \ref{Corollary 3.4}, the function of $t$
\begin{equation*}
\||x-y|e^{it\triangle}u_0(x)\,e^{it\triangle}v_0(y)\|_{L^2(\mathbb R^{2n}_{x,y})}^2
\end{equation*}
is convex in $\mathbb R$.
\end{corollary}
\begin{proof} For $\gamma >0$, choose $\alpha =0$ and $\beta=\sqrt{\frac 2\gamma}$ in \eqref{dos3}. It shows that
\begin{equation*}
\|e^{\frac\gamma 2 |x-y|^2}e^{it\triangle}u_0(x)\,e^{it\triangle}v_0(y)\|_{L^2(\mathbb R^{2n}_{x,y})}^2
\end{equation*}
is logarithmically convex in $[0,T]$. This and the fact that $$
\|e^{it\triangle}u_0(x)\,e^{it\triangle}v_0(y)\|_{L^2(\mathbb R^{2n}_{x,y})},
$$ 
is constant, given that
\begin{multline}\label{E: antesdellimite}
\int_{\mathbb R^{2n}}\frac{e^{\gamma |x-y|^2}-1}{\gamma}|e^{it\triangle}u_0(x)\,e^{it\triangle}v_0(y)|^2\,dxdy\\\le \left(1-\frac tT\right)\int_{\mathbb R^{2n}}\frac{e^{\gamma |x-y|^2}-1}{\gamma}|u_0(x)\,v_0(y)|^2\,dxdy\\+\frac tT\int_{\mathbb R^{2n}}\frac{e^{\gamma |x-y|^2}-1}{\gamma}|e^{iT\triangle}u_0(x)\,e^{iT\triangle}v_0(y)|^2\,dxdy,
\end{multline}
 when $0\le t\le T$. The corollary follows from \eqref{E: antesdellimite}, after letting $\gamma$ tend to zero. 
\end{proof}
\begin{remark}\label{quizas} 
Corollary \ref{Corollary Viriel interactive identity} follows from \eqref{dos3} because the constant in front of the right hand side of the inequality is precisely equal to $1$.

Similar arguments show  that
\begin{equation*}
\||x-y|^{\frac\alpha 2}e^{it\triangle}u_0(x)\,e^{it\triangle}v_0(y)\|_{L^2(\mathbb R^{2n}_{x,y})}^2
\end{equation*}
is convex, when $n\ge 3$ and $1\le \alpha\le 2$. In \cite{CoKeStTaTa} the authors proved a similar result, when $u_0=v_0$, also in the non-linear setting.
\end{remark}

\subsection{Galilean invariance}  
  
  The following result, which is a consequence of Theorem \ref{Theorem 1.2} and the Galilean invariant property of the  free Schr\"odinger group, describes  the time evolution  of the location of the \lq\lq mass"  of a Gaussian decaying solution.

\begin{corollary}
 \label{Corollary 3.5} 
 
 Using the same hypotheses and notation that in Theorem \ref{Theorem 1.2}. For any $\,\nu\in\mathbb R^n$
 \begin{equation}
 \label{gal1}
\|e^{\frac{|x+2t\nu|^2}{(\alpha t+\beta)^2}} u(t)\|_2
\leq \|e^{\frac{|x|^2}{\beta^2}}u(0)\|_2^{\theta(t)}
\|e^{\frac{|x+2T\nu|^2}{(\alpha T+\beta)^2}} u(T)\|_2^{1-\theta(t)},
\end{equation}
when $0\le t\le T$, and
\begin{equation}
\label{gal2}
\|e^{\frac{|x+t\nu|^2}{(\alpha t+\beta)^2}} u(t)\|_2
\leq \|e^{\frac{|x|^2}{\beta^2}}u(0)\|_2^{\mu(t)}
\|e^{\frac{4|\xi+\nu|^2}{\alpha^2}}\widehat u(0)\|_2^{1-\mu(t)},
\end{equation}
when $t\ge 0$.
\end{corollary}

\begin{proof} We recall the Galilean invariance  of the free Schr\"odinger group
\begin{equation}
\begin{aligned}
\label{gal3}
u_{\nu}(x,t)=&e^{it\Delta}(e^{i\nu\cdot \,}u_0(\cdot))(x)=e^{-i|\nu|^2 t} e^{i\nu\cdot x}(e^{it\Delta}u_0)(x-2t\nu)\\
=&e^{-i|\nu|^2 t} e^{i\nu\cdot x} u_0(x-2t\nu,t).
\end{aligned}
\end{equation}

Thus combining \eqref{gal3}, the identity
$$
\widehat{ e^{i\nu\cdot x\,}u_0}(\xi)=u_0(\xi-\nu),
$$
and the inequalities \eqref{2.22}-\eqref{2.23} in Theorem \ref{Theorem 1.2}, we obtain \eqref{gal1} and \eqref{gal2}.
\end{proof}

\subsection{Final remarks}
Next, we recall the  following result established in \cite{KePoVe}, which is one of the main estimates in that paper :
\begin{lemma}\label{L: ultimo lemma} There exists $\epsilon>0$ such that if 
\begin{equation}
\label{hyp2}
 V:\mathbb R^n\times [0,T]\to\mathbb C,\;\;\;\;\text{with}\;\;\;\;
\|V\|_{L^1_tL^{\infty}_x}\leq \epsilon,
\end{equation}
and $u\in C([0,T]:L^2_x(\mathbb R^n))$ is a strong solution of the IVP 
\begin{equation}
\begin{cases}
\begin{aligned}
\label{eq1}
&\p_tu=i(\Delta +V(x,t))u+F(x,t),\\
&u(x,0)=u_0(x).
\end{aligned}
\end{cases}
\end{equation}
with
\begin{equation}
\label{hyp3} u_0,\,u_T\equiv u(\,\cdot\,,T)\in
L^2(e^{2\lambda\cdot x}dx),\;F\in L^1([0,T]:L^2_x(e^{2\lambda\cdot
x}dx)),
\end{equation}
for some $\lambda\in\mathbb R$, then there exists $c$ independent of 
$\lambda$ such that
\begin{equation}
\label{uno}
\sup_{0\leq t\leq T}\| e^{\lambda\cdot x} u(\,\cdot\,,t)\|_2 \leq c
\Big(\|e^{\lambda\cdot x} u_0\|_2 + \|e^{\lambda\cdot x} u_T\|_2 +\int_0^T
\|e^{\lambda\cdot x} F(\cdot, t)\|_2 dt\Big).\;\;\square
\end{equation}

\end{lemma}

Notice that in the above result one assumes the existence of a reference
$L^2$-solution $u$ of the equation \eqref{eq1} and then, under the 
hypotheses \eqref{hyp2} and \eqref{hyp3}, shows that the exponential decay in the
time interval $[0,T]$ is preserved. From the arguments used in \eqref{key1}-\eqref{explain2}, it follows
that the inequality \eqref{uno} holds with Gaussian weights, i.e. with $\gamma |x|^2$, $\gamma >0$ instead of $\lambda\cdot x$ in the exponent. 
 
  In a forthcoming work on Schr\"odinger equations with variable coefficients, among other results, we shall extend those in Lemma \ref{L: ultimo lemma} to a class of potentials $V$ without smallness assumptions, and give suitable density arguments to justify the manipulations, integrations by parts and calculations, which arise at the time of trying to derive, in this more general context, both the preservation and the logarithmic convexity  of the $L^2(e^{2\gamma |x|^2}dx)$-norm, $\gamma >0$,  of the corresponding solutions.


\begin{thebibliography}{10} 


\bibitem{CoKeStTaTa} Colliander, J., Keel, M., Staffilani, G.,  Takaoka, H., and Tao, T.,
{Global well-posedness and scattering in the energy space for the critical nonlinear Schr\"odinger equation in $R^3$}. To appear in Annals of Math.


 \bibitem{CoPr} Cowling, M., and Price, J. F., 
{Generalizations of Heisenberg's inequality}, Harmonic Analysis (Cortona, 1982) Lecture Notes in Math.,{\bf 992} (1983), 443-449, Springer, Berlin

 
 
\bibitem{EsKePoVe}  Escauriaza, L., Kenig, C. E., Ponce, G., and Vega, L., {On Uniqueness Properties of Solutions of Schr\"odinger  Equations}, Comm. PDE. {\bf 31} (2006), 1811-1823

\bibitem{HaNaTs}  Hayashi, N., Nakamitsu, K., and Tsutsumi, N., {On solutions of the initial value problem  for the nonlinear Schr\"odinger equations}, J. Funct. Anal., {\bf 71} (1987), 218-245

\bibitem{Ho} H\"ormander, L., {A uniqueness theorem of Beurling for Fourier transform pairs}, 
Ark. Mat. 29 (1991), no. 2, 237--240

\bibitem{IoKe} Ionescu, A. D., and Kenig, C. E., {$L^p$ Carleman inequalities and uniqueness of solutions of nonlinear Schr\"odinger equations}, Acta Math. {\bf193} (2004), 193-239



\bibitem{KePoVe} Kenig, C. E., Ponce, G., and Vega, L.,
{On the support of solutions of nonlinear Schr\"odinger equations}, Comm. Pure Appl. Math. {\bf 60}
(2002), 1247-1262

\bibitem{SiSu} Sitaram, A., Sundari, M., and Thangavelu, S., {Uncertainty principles on certain Lie groups},  Proc. Indian Acad. Sci. Math. Sci. {\bf 105}
(1995), 135-151


\bibitem{StSh} Stein, E. M., and Shakarchi. R., {Complex Analysis},  Princeton Lectures in Analysis, Princeton University Press,  (2003)




\end{thebibliography}
\end{document}